\documentclass[reqno]{amsart}

\usepackage{amscd}

\theoremstyle{plain}
  \newtheorem{theo}{Theorem}
  
  \newtheorem{prop}{Proposition}
  \newtheorem{lemm}{Lemma}
  \newtheorem*{quest*}{Question}  

\theoremstyle{remark}
  \newtheorem*{rema*}{Remark}
  \newtheorem{rema}{Remark}
  \newtheorem{exam}{Example}

\newcommand\itemref[1]{(\ref{#1})}

\newcommand\ie{i.e.\ }
\newcommand\cf{cf.\ }

\renewcommand\o{\circ}
\newcommand\goe{\mathfrak g}

\newcommand\aut{\text{\rm aut}}
\newcommand\ham{\text{\rm ham}}
\newcommand\Aut{\text{\rm Aut}}
\newcommand\Emb{\text{\rm Emb}}
\newcommand\Diff{\text{\rm Diff}}

\newcommand\id{\text{\rm id}}

\newcommand\vol{\nu}
\newcommand\OGr{\text{\rm Gr}}

\newcommand\SOGr{\text{\rm S}\OGr}
\newcommand\Ham{\text{\rm Ham}}

\DeclareMathOperator\tr{tr}
\DeclareMathOperator\grad{grad}
\DeclareMathOperator\II{II}
 
\newcommand\vf{\mathfrak X}

\newcommand\Diffvol{\text{\rm Diff}(M,\vol)}

\newcommand\Hamvol{\Ham(M,\vol)}
\newcommand\hamvol{\ham(M,\vol)}

\newcommand\R{\mathbb R}
\newcommand\Z{\mathbb Z}
\newcommand\N{\mathbb N}

\begin{document}

\title{Non-linear {G}rassmannians as coadjoint orbits}

\author{Stefan Haller}

\address{Stefan Haller, 
         Department of Mathematics, University of Vienna,
         Strudlhofgasse 4, A-1090 Vienna, Austria.}

\email{Stefan.Haller@univie.ac.at}

\author{Cornelia Vizman}

\address{Cornelia Vizman, 
         West University of Timisoara, Department of Mathematics, 
         Bd. V.Parvan 4, 1900 Timisoara, Romania.}

\email{vizman@math.uvt.ro}

\thanks{Both authors are supported by the `Fonds zur F\"orderung der
        wissenschaftlichen Forschung' (Austrian Science Fund),
        project number {\tt P14195-MAT}}



\begin{abstract}
For a given manifold $M$ we consider the non-linear Grassmann manifold
$\OGr_n(M)$ of $n$--dimensional submanifolds in $M$.
A closed $(n+2)$--form on $M$ gives rise to a closed $2$--form on
$\OGr_n(M)$. If the original form was integral, the $2$--form will be the
curvature of a principal $S^1$--bundle over $\OGr_n(M)$. Using this 
$S^1$--bundle one obtains central extensions for certain groups of 
diffeomorphisms of $M$. We can realize $\OGr_{m-2}(M)$ 
as coadjoint orbits of the extended group of exact volume
preserving diffeomorphisms and the symplectic Grassmannians 
$\SOGr_{2k}(M)$ as coadjoint
orbits in the group of Hamiltonian diffeomorphisms.
We also generalize the vortex filament equation 
as a Hamiltonian equation on $\OGr_{m-2}(M)$. 
\end{abstract}

\maketitle

\section{Introduction}

Let $M$ be a smooth connected closed manifold of dimension $m$.
We are interested in the
space of closed submanifolds of $M$. More precisely we fix a dimension $n$
and let $\OGr_n(M)$ denote the space of all $n$--dimensional oriented compact
boundaryless submanifolds of $M$. This is easily seen to be a Fr\'echet 
manifold in a natural way. We consider this as a non-linear analogue of the 
classical Grassmann manifolds. 

Every closed differential form $\alpha$ of degree $n+2$ on $M$ gives rise to 
a closed $2$--form $\tilde\alpha$ on $\OGr_n(M)$. If $\alpha$ was 
integral, our first theorem says that there is a principal $S^1$--bundle
$\mathcal P\to\OGr_n(M)$ with principal connection, whose curvature form 
is $\tilde\alpha$. The group of equivariant connection preserving
diffeomorphisms of $\mathcal P$ then is a central extension of the group of
Hamiltonian diffeomorphisms on $\OGr_n(M)$. The latter makes sense, even if
$\tilde\alpha$ is degenerate. Restricting everything to a connected component
of $\OGr_n(M)$, the extension becomes $1$--dimensional with fiber $S^1$.

Now the group of diffeomorphisms of $M$ which preserve $\alpha$ acts 
symplectically on $\OGr_n(M)$. In some cases there are interesting
subgroups $G$, which actually act in a Hamiltonian way. In such a 
situation the pull back of the central extension described above gives
a central extension $1\to S^1\to\tilde G\to G\to1$. 

Let us describe two cases in more detail. Suppose $\alpha$ was an integral volume
form and let $\OGr_{m-2}(M)$ denote the space of codimension $2$ submanifolds.
As mentioned above, the volume form gives a closed $2$--form on $\OGr_{m-2}(M)$,
which turns out to be (weakly) non-degenerate. Then the group of exact
volume preserving diffeomorphisms acts in a Hamiltonian way on $\OGr_{m-2}(M)$.
So the pull back gives central extensions $\tilde G$ of the group of exact 
volume preserving diffeomorphisms by $S^1$. This is Ismagilov's way of constructing
these extensions, see \cite{I96}. Using the moment map we will then realize
the symplectic manifold $\OGr_{m-2}(M)$ as a coadjoint orbit of the group
$\tilde G$. We even get Lie group structure on the extensions
$\tilde G$ of the group of exact volume preserving diffeomorphisms.

For the second situation we have in mind we start with a symplectic 
manifold $(M,\omega)$. Taking $\alpha:=\omega^{k+1}$ we get a closed $2$--form
$\tilde\alpha$ on $\OGr_{2k}(M)$. This form is no longer symplectic.
However, when restricted to the open subset $\SOGr_{2k}(M)$ of symplectic 
submanifolds it will become non-degenerate. We will refer to $\SOGr_{2k}(M)$
as a non-linear symplectic Grassmannian. The group of Hamiltonian 
diffeomorphisms of $M$ now acts in a Hamiltonian way on $\SOGr_{2k}(M)$. So 
the procedure above yields central extensions of $\Ham(M,\omega)$ by $S^1$.
These extensions are not very interesting, since the associated
extensions of Lie algebras turn out to be trivial. However, it permits us to 
realize the symplectic manifold $\SOGr_{2k}(M)$ as a coadjoint orbit of
$\Ham(M,\omega)$.

For a Riemannian manifold $(M,g)$ the non-linear Grassmannian $\OGr_{m-2}(M)$ 
of codimension two submanifolds 
has a canonical almost K\"ahler structure.
The $g$--volume of the submanifold gives a smooth function on $\OGr_{m-2}(M)$
and its Hamiltonian equation generalizes the vortex filament equation.

Finally, let us remark that everything generalizes to non-compact $M$ in a
straight forward way. The diffeomorphism groups then have to be replaced by
the compactly supported ones.

\section{Non-linear Grassmannians}

Throughout the whole paper $M$ will be a 
smooth closed connected $m$--dimensional manifold. 
Let $\OGr_n(M)$ denote the space of all oriented compact  
$n$--dimensional not necessarily connected submanifolds without boundary. 
This is easily seen to be a Fr\'echet manifold in a natural way, see 
\cite{KM97}.
Note that there is a natural action of the group $\Diff(M)$ on 
$\OGr_n(M)$. A classical theorem due to R.~Thom implies that $\Diff(M)_0$,
the connected component in the group of diffeomorphisms,
acts transitively on every connected component of $\OGr_n(M)$.

Suppose $N\in\OGr_n(M)$. Then the tangent space of $\OGr_n(M)$ at $N$ can
naturally be identified with the space of smooth sections of the normal
bundle $TN^\perp:=TM|_N/TN$. Any $\alpha\in\Omega^k(M)$ gives rise
to $\tilde\alpha\in\Omega^{k-n}(\OGr_n(M))$ via:
$$
(\tilde\alpha)_N(Y_1,\dotsc,Y_{k-n})
:=\int_Ni_{Y_{k-n}}\cdots i_{Y_1}\alpha.
$$
Here $N\in\OGr_n(M)$ and $Y_j$ are tangent vectors at $N$, \ie sections
of $TN^\perp$. Then $i_{Y_{k-n}}\cdots i_{Y_1}\alpha\in\Omega^n(N)$
does not depend on representatives $Y_j$ and integration is well defined, for 
$N\in\OGr_n(M)$ comes with an orientation. 

Let $\zeta$ denote the infinitesimal $\Diff(M)$--action on $\OGr_n(M)$, that
is for every vector field $X\in\vf(M)$ on $M$ we have a fundamental vector
field $\zeta_X$ on $\OGr_n(M)$.
One easily verifies the following 

\begin{lemm}
For every $X\in\vf(M)$, $N\in\OGr_n(M)$, $k\in\N$, $\alpha\in\Omega^k(M)$
and every $\varphi\in\Diff(M)$ we have:
\begin{enumerate}
\item
$\zeta_X(N)=X|_N$.
\item
$\widetilde{d\alpha}=d\tilde\alpha$.
\item
$i_{\zeta_X}\tilde\alpha=\widetilde{i_X\alpha}$.
\item
$L_{\zeta_X}\tilde\alpha=\widetilde{L_X\alpha}$.
\item
$\varphi^*\tilde\alpha=\widetilde{\varphi^*\alpha}$.
\end{enumerate}
\end{lemm}

Suppose $\alpha$ is a closed $k$--form on $M$. Then we get a closed
$2$--form $\tilde\alpha$ on $\OGr_{k-2}(M)$. Our first theorem states that
if $[\alpha]\in H^k(M;\R)$ is integral then $\tilde\alpha$ will be the
curvature form of a principal $S^1$--bundle over
$\OGr_{k-2}(M)$.

\begin{theo}\label{main-theo}
Let $M$ be a closed connected manifold and let $\alpha$ be a closed 
$k$--form representing an integral cohomology class of $M$.
Then there exist a principal $S^1$--bundle
$\mathcal P\to\OGr_{k-2}(M)$ and a principal
connection $\eta\in\Omega^1(\mathcal P)$ whose curvature form
is $\tilde\alpha$. 
\end{theo}

\begin{proof}
Note that it suffices to proof this for one representative of
$[\alpha]\in H^k(M;\R)$. For if $\eta$ is a principal connection
with curvature $\tilde\alpha$ then $\eta+\pi^*\tilde\beta$
is a principal connection with curvature form $\widetilde{\alpha+d\beta}$.

Pick a smooth triangulation $\Delta^*$ of $M$ and let $\Delta^{m-k}$
denote its $(m-k)$--skeleton. Choose an open neighborhood $U$ of
$\Delta^{m-k}$ which deformation retracts onto $\Delta^{m-k}$.
Moreover set $A:=M\setminus U$ and $V:=A^\circ=M\setminus\overline U$.
One easily checks the following properties:
\begin{enumerate}
\item\label{prop-i}
$H^j(A;\Z)=0$ for all $j\geq k$. We will actually  only use $H^k(A;\Z)=0$.
\item\label{prop-ii}
For all compact $K\subseteq M\setminus\Delta^{m-k}$ there exists
$g\in\Diff(M)_0$ with $g(K)\subseteq V$. If
moreover $K'\subseteq V$ compact, then $g$ and the diffeotopy connecting
it to the identity can be chosen to fix the points in $K'$.
\end{enumerate}
Using \itemref{prop-i} and considering
$$
\CD
H^k(M,A;\R) & @>>> & H^k(M;\R) &
\\
@AAA        &      & @AAA &
\\
H^k(M,A;\Z) & @>>> & H^k(M;\Z) & @>>> H^k(A;\Z)
\endCD
$$
we see that $[\alpha]\in H^k(M;\R)$ has a representative
which vanishes on $V$ and which represents an integral class
in $H^k(M,A;\R)$, \ie lies in the image of $H^k(M,A;\Z)\to H^k(M,A;\R)$.
Since it suffices to construct the bundle and the
connection for some representative we may assume from now on
\begin{enumerate}
\setcounter{enumi}{2}
\item\label{vanishes_on_V}
$\alpha$ vanishes on $V$.
\item\label{rel_integral} 
$\alpha$ represents an integral class in $H^k(M,A;\R)$.
\end{enumerate}

\begin{lemm}\label{easy-lemma}
Suppose $L$ is a compact manifold of dimension $l<k$, which might have a
boundary (even corners) and suppose $f:L\to M$ smooth. Then there exists
$g\in\Diff(M)_0$, such that $g(f(L))\subseteq V$.
Moreover if $f(\partial L)\subseteq V$ then $g$ and the diffeotopy
connecting it with the identity can be chosen to fix the points in 
$f(\partial L)$.
\end{lemm}

\begin{proof}[Proof of Lemma \ref{easy-lemma}]
A well known transversality argument shows that there
exists $g_1\in\Diff(M)_0$ with $g_1\circ f$ transversal to
$\Delta^{m-k}$. Since $l<k$ we thus must have
$g_1(f(L))\cap\Delta^{m-k}=\emptyset$. From \itemref{prop-ii} we get 
$g_2\in\Diff(M)_0$ with $g_2(g_1(f(L))\subseteq V$.
The second part is proved similarly.
\end{proof}

Let us continue with the proof of Theorem~\ref{main-theo}.
Let $I=[0,1]$ denote the unit interval.
For $\varphi\in C^\infty(I,\Diff(M))$ with $\varphi_0=\id$ we
define 
$$
\mathcal U_\varphi:=\big\{N\in\OGr_{k-2}(M):\varphi_1(N)\subseteq V\big\}.
$$
Lemma~\ref{easy-lemma} implies that 
$\mathcal U_\varphi$ is an open
covering of $\OGr_{k-2}(M)$. Moreover we set
$$
\lambda_\varphi
:=-\int_0^1\varphi_t^*i_{\dot\varphi_t}\alpha\ dt
\in\Omega^{k-1}(M).
$$
A one line computation shows
$d\lambda_\varphi=\alpha-\varphi^*_1\alpha$. Because of
\itemref{vanishes_on_V} we particularly have
$d\lambda_\varphi=\alpha$ on $\varphi^{-1}_1(V)$ and thus
$d\widetilde{\lambda_\varphi}=\tilde\alpha$ on $\mathcal U_\varphi$.
The $\widetilde{\lambda_\varphi}$ will be the connection forms, we are going to define
the transition cocycle defining the $S^1$--bundle.

Fix $\varphi,\psi\in C^\infty(I,\Diff(M))$ with $\varphi_0=\psi_0=\id$ and
consider homotopies $\Phi\in C^\infty(I\times I,\Diff(M))$
with $\Phi_{0,t}=\varphi_t$, $\Phi_{1,t}=\psi_t$ and
$\Phi_{s,0}=\id$ for all $s\in I$. For such a $\Phi$ we set
$$
\mathcal U_\Phi:=\big\{N\in\OGr_{k-2}(M):
\text{$\Phi_{s,1}(N)\subseteq V$ for all $s\in I$}\big\}.
$$
clearly $\mathcal U_\Phi$ are open subsets of $\mathcal U_\varphi\cap\mathcal U_\psi$.
Lemma~\ref{easy-lemma} shows that for every $N\in\mathcal U_\varphi\cap\mathcal U_\psi$ there 
exists a homotopy $\Phi$ with ends $\varphi$ and $\psi$, such that 
$N\in\mathcal U_\Phi$. In other words $\mathcal U_\Phi$ constitute an open covering 
of $\mathcal U_\varphi\cap\mathcal U_\psi$, as $\Phi$ varies with fixed ends
$\varphi$ and $\psi$. For such a $\Phi$ we define
$$
\tau_\Phi:=\int_0^1\int_0^1\Phi_{s,t}^*i_{\delta\Phi(\partial_t)}
i_{\delta\Phi(\partial_s)}\alpha\ ds\ dt\in\Omega^{k-2}(M).
$$
Using the Maurer--Cartan equation for the left logarithmic derivative,
\cf~\cite{KM97},
$$
[\delta\Phi(\partial_s),\delta\Phi(\partial_t)]=
\frac\partial{\partial t}\delta\Phi(\partial_s)
-\frac\partial{\partial s}\delta\Phi(\partial_t),
$$
an easy computation yields
$$
d\tau_\Phi=\lambda_\psi-\lambda_\varphi
+\int_0^1\Phi^*_{s,1}i_{\delta\Phi(\partial_s)}\alpha\ ds.
$$
Particularly $\widetilde{\lambda_\psi}-\widetilde{\lambda_\varphi}
=d\widetilde{\tau_\Phi}$ on $\mathcal U_\Phi$. 
Note that $\widetilde{\tau_\Phi}$ is a function on $\OGr_{k-2}(M)$ and obviously
$\widetilde{\tau_\Phi}(N)=\int_{I\times I\times N}\hat\Phi_{N}^*\alpha$
with $\hat\Phi_{N}(s,t,x)=\Phi_{s,t}(x)$, $x\in N$. 
If $\Psi$ is another homotopy with
ends $\varphi$ and $\psi$ and $N\in\mathcal U_\Phi\cap\mathcal U_\Psi$ then
$$
\widetilde{\tau_\Psi}(N)-\widetilde{\tau_\Phi}(N)\in\Z.
$$
Indeed, $\hat\Psi_{N}-\hat\Phi_{N}$ represents a class in 
$H_k(M,A;\Z)$ and $\widetilde{\tau_\Psi}(N)-\widetilde{\tau_\Phi}(N)$ 
is the pairing of this class with $[\alpha]\in H^k(M,A;\R)$. 
From \itemref{rel_integral} we see that
the result must be integral as well.

So when considered as functions $\mathcal U_\Phi\to S^1:=\R/\Z$ the
$\widetilde{\tau_\Phi}$ fit together and define well defined smooth 
$f_{\varphi\psi}:\mathcal U_\varphi\cap\mathcal U_\psi\to S^1$
satisfying 
$df_{\varphi\psi}=\widetilde{\lambda_\psi}-\widetilde{\lambda_\varphi}$.
A similar argument shows that they satisfy the cocycle condition
$f_{\varphi\psi}+f_{\psi\rho}-f_{\varphi\rho}=0$ as functions
$\mathcal U_\varphi\cap\mathcal U_\psi\cap\mathcal U_\rho\to S^1$, where $S^1$
is written additively.

Now define $\mathcal P$ to be the principal $S^1$--bundle one obtains when
gluing $\mathcal U_\varphi\times S^1$ with the help of $f_{\varphi\psi}$.
On $\mathcal U_\varphi\times S^1$ we
define
$\eta_\varphi:=\lambda_\varphi+d\theta$, where $d\theta$ denotes the
standard volume form on $S^1$. These locally defined $\eta_\varphi$
define a global principal connection $\eta\in\Omega^1(\mathcal P)$, for we
have $df_{\varphi\psi}=\widetilde{\lambda_\psi}
-\widetilde{\lambda_\varphi}$ on $\mathcal U_\varphi\cap\mathcal U_\psi$. Since
$d\widetilde{\lambda_\varphi}=\tilde\alpha$ on $\mathcal U_\varphi$ its 
curvature form is $\tilde\alpha$.
This finishes the proof of the theorem.
\end{proof}


\begin{exam}\label{trivial_ex}
Let us consider the case $k=2$. So $\alpha$ is a closed integral
$2$--form and $\OGr_{k-2}(M)$ is the space of oriented points in $M$. 
Let $\mathcal M$ denote the connected component of $\OGr_{k-2}(M)$
where the submanifolds consist of a single positively oriented point.
Certainly $\mathcal M=M$ and $\tilde\alpha=\alpha$.
So in this case the restriction of the bundle $\mathcal P\to\OGr_{k-2}(M)$ to 
$\mathcal M$ gives the classical circle bundle with connection 
corresponding to the closed integral $2$--form $\alpha$.
\end{exam}

\begin{rema}
A Theorem of R.~Thom implies that the action of $\Diff(M)_0$ on connected 
components of $\OGr_n(M)$ is transitive. Hence connected components of 
$\OGr_n(M)$ are homogeneous spaces of $\Diff(M)_0$.
Actually Thom's theorem shows that $\Diff(M)_0$ acts transitively on
connected components of 
$\Emb(N,M)$, the space of smooth embeddings of $N$ in $M$.
So connected components of $\Emb(N,M)$ are homogeneous spaces of
$\Diff(M)_0$ too.
Below we will see that similar statements hold for the group of 
volume preserving diffeomorphisms.

Moreover the connected components of $\Emb(N,M)$ are principal
bundles over corresponding connected components of
$\OGr_n(M)$, see \cite{KM97}. The structure group is 
the group of orientation preserving diffeomorphisms of $N$.
\end{rema}

\section{Universal construction}\label{univ_constr}

Suppose we have a principal $S^1$--bundle
$\pi:(\mathcal P,\eta)\to(\mathcal M,\Omega)$ with connection $\eta$ and
curvature $\Omega$. We assume $\mathcal M$ connected but it may very well be
infinite dimensional.
We associate Kostant's exact sequence of groups, see \cite{Ko70}:
$$
1\to S^1\to\Aut(\mathcal P,\eta)\to\Ham(\mathcal M,\Omega)\to 1.
$$
Here $\Aut(\mathcal P,\eta)$ is the connected component of
the group of equivariant connection preserving 
diffeomorphisms of $\mathcal P$ and $\Ham(\mathcal M,\Omega)$ is
the group of Hamiltonian diffeomorphisms of $\mathcal M$. The latter can
either be described as the connected component of holonomy preserving 
diffeomorphisms, or as the kernel of a flux homomorphism  \cite{NV03}.

The group $\Aut(\mathcal P,\eta)$ acts on $\mathcal M$ in a Hamiltonian way
with equivariant moment map
$$
\hat\mu:\mathcal M\to\aut(\mathcal P,\eta)^*,
\qquad
\hat\mu(x)(\xi)=-(i_\xi\eta)(\pi^{-1}(x)).
$$
This moment map is universal in the following sense:
Whenever we have a Hamiltonian action of a Lie group $G$ on $\mathcal M$,
we can pull back Kostant's extension and get a $1$--dimensional central 
group extension:
$$
\CD
S^1 & @>>> & \Aut(\mathcal P,\eta) & @>>> & \Ham(\mathcal M,\Omega)&
\\
@|         &      & @AAA & & @AAA &
\\
S^1 & @>>> & \tilde G & @>>> & G
\endCD
$$
This is a Lie group extension, even if Kostant's extension
is only a group extension in this infinite dimensional setting, 
see \cite {NV03}.
Moreover the pull back $\tilde\mu:\mathcal M\to\tilde\goe^*$ of 
$\hat\mu$ is a smooth equivariant
moment map for the $\tilde G$--action on $\mathcal M$.
Consider the corresponding central extension of Lie algebras:
$$
\CD
\R & @>>> & \aut(\mathcal P,\eta) & @>>> & \ham(\mathcal M,\Omega)
\\
@|        &      & @AAA & & @AA\zeta A &
\\
\R & @>>> & \tilde\goe & @>p>> & \goe
\endCD
$$

\begin{prop}\label{kks}
In the situation above suppose moreover that $G$ acts transitively on 
$\mathcal M$ and admits an injective but not necessarily equivariant moment 
map $\mu:\mathcal M\to\goe^*$.
Then the equivariant moment map $\tilde\mu:\mathcal M\to\tilde\goe^*$
is one-to-one onto a coadjoint orbit of $\tilde G$. Moreover it pulls back
the Kirillov--Kostant--Souriau symplectic form to $\Omega$.
\end{prop}

\begin{proof}
Note first, that $p^*\circ\mu:\mathcal M\to\tilde\goe^*$ is an injective 
but not necessarily equivariant moment map for the $\tilde G$--action on 
$\mathcal M$. Since $\mathcal M$ is connected, two moment maps differ 
by a constant in $\tilde\goe^*$. Thus every moment map for the $\tilde
G$--action on $\mathcal M$ is injective, particularly $\tilde\mu$.
Next, $\tilde G$ acts transitively, for $G$ does. So the equivariance of
$\tilde\mu$ implies that $\tilde\mu$ is onto a single coadjoint orbit.
A straight forward calculation shows that the pull back of the
Kirillov--Kostant--Souriau symplectic form is $\Omega$.
\end{proof}

Till the end of the section we will denote all the left $G$--actions
by a dot.
Suppose we have a not necessarily equivariant moment map 
$\mu:\mathcal M\to\goe^*$. Let $h:\goe\to C^\infty(\mathcal M,\R)$ denote
the dual map, that is $h_X(x)=\mu(x)(X)$, for $x\in\mathcal M$ and
$X\in\goe$. The universal property of the pull back implies that there is a
unique section $\sigma:\goe\to\tilde\goe$ with $i_{\sigma(X)}\eta=-\pi^*h_X$.
Conversely every section is obtained in this way. So we have a one-to-one
correspondence of not necessarily equivariant moment maps $\mu:\mathcal
M\to\goe^*$ and sections of $p:\tilde\goe\to\goe$. Every such choice gives
a linear isomorphism 
\begin{equation}\label{split}
\R\oplus\goe\to\tilde\goe,
\quad
(a,X)\mapsto a+\sigma(X).
\end{equation}
Via \eqref{split} the equivariant moment map 
$\tilde\mu:\mathcal M\to\tilde\goe^*$ we constructed above is
$$
\tilde\mu:\mathcal M\to(\R\oplus\goe)^*=\R^*\oplus\goe^*,
\quad
\tilde\mu=(-1^*,\mu).
$$
Here $1^*$ is the dual base to $1$ considered as base of $\R$.
Equivalently $\tilde\mu(x)(a,X)=\mu(x)(X)-a$, for $x\in\mathcal M$,
$X\in\goe$ and $a\in\R$.

Define $\kappa:G\to\goe^*$ by 
$-\kappa(g^{-1})(X)=g\cdot\sigma(X)-\sigma(g\cdot X)$. 
So $\kappa$ is the failure of
$\sigma$ to be $G$--equivariant. Then via \eqref{split} the adjoint action 
is
$$
g\cdot(a,X)=(a-\kappa(g^{-1})(X),g\cdot X),
\quad\text{for $g\in G$.}
$$
The function $\kappa$ satisfies
$\kappa(g_1g_2)=\kappa(g_1)+g_1\cdot \kappa(g_2)$,
hence it is a 1--cocycle (derivation) on $G$ with values in $\goe^*$.
For $g\in G$ and $X\in\goe$ the function
$h_{g\cdot X}-g\cdot h_X$ is locally constant, 
hence constant since $\mathcal M$ is connected. 
So we get a
function $G\to\goe^*$ which measures the failure of the moment map to be
$G$--equivariant. One readily checks  
$-\kappa(g^{-1})(X)=h_{g\cdot X}-g\cdot h_X$, 
equivalently $\kappa(g)=\mu(g\cdot x_0)-g\cdot\mu(x_0)$,
for every $x_0\in\mathcal M$.
So the section $\sigma$ is $G$--equivariant iff
the corresponding moment map is $G$--equivariant.

Via $\eqref{split}$ we can express the Lie bracket as
$$
\bigl[(a,X),(b,Y)\bigr]=\bigl(c(X,Y),[X,Y]\bigr),
$$
where $c\in\Lambda^2\goe^*$ is the cocycle 
$c(X,Y)=[\sigma(X),\sigma(Y)]-\sigma([X,Y])$. Note that $c$ also is a
measure for the failure of $\sigma$ to be $\goe$--equivariant.
By choosing different sections $\sigma$ we obtain all 2--cocycles $c$
in one cohomology class, but different sections could define
the same 2--cocycle. 
Moreover
the differential of $\kappa:G\to\goe^*$ at the identity satisfies
$(T_e\kappa\cdot X)(Y)=c(X,Y)$. Since we had 
$-\kappa(g^{-1})(X)=h_{g\cdot X}-g\cdot h_X$, we get
\begin{equation}\label{cocycle}
c(X,Y)=h_{[X,Y]}+L_{\zeta_X}h_Y=h_{[X,Y]}+\{h_X,h_Y\}
=h_{[X,Y]}-\Omega(\zeta_X,\zeta_Y).
\end{equation}
The unexpected signs of the second summands stem from the convention for 
the Lie derivative of functions, which is an infinitesimal right action
and quite confusing.
Thus $c$ also is a measure for the failure
of the moment map to be $\goe$--equivariant. Particularly the moment map
$\mu$ is $\goe$--equivariant iff the corresponding section $\sigma$ is
$\goe$--equivariant.
Finally for every point $x_0\in\mathcal M$ we have 
$c(X,Y)=h_{[X,Y]}(x_0)-\Omega(\zeta_X,\zeta_Y)(x_0)$.
So we see that $c(X,Y)=-\Omega(\zeta_X,\zeta_Y)(x_0)$ is a cocycle describing
the extension $0\to\R\to\tilde\goe\to\goe\to0$ and corresponds to moment 
maps satisfying $\mu(x_0)=0$.

\section{Codimension two Grassmannians}

Let $M$ be a closed $m$--dimensional manifold with integral volume form 
$\vol$, that is $\int_M\vol\in\Z$. From Theorem~\ref{main-theo} we get 
a principal $S^1$--bundle
$\mathcal P\to\OGr_{m-2}(M)$ and a principal
connection $\eta$ whose curvature form
is $\Omega:=\tilde\vol$. 
Recall that $\Omega_N( Y _1, Y _2)=\int_Ni_{ Y _2}i_{ Y _1}\vol$
for tangent vectors $ Y _1$ and $ Y _2$ at $N$, \ie sections of $TN^\perp$. 
Note that $\Omega$ is symplectic, \ie (weakly) non-degenerate.
The action of the group of volume preserving
diffeomorphisms $\Diffvol$ on $\OGr_{m-2}(M)$ 
preserves the symplectic form $\Omega$.
In dimension $m=3$ the symplectic form $\Omega$ is known as the 
Marsden--Weinstein 
symplectic from on the space of unparameterized oriented links, see
\cite{MW83}.

Let $\Ham(M,\vol)$ denote the group of exact volume preserving
diffeomorphisms with Lie algebra
$$
\ham(M,\vol)
=\{X\in\vf(M):\text{$i_X\vol$ exact differential form}\}.
$$
The action of $\Ham(M,\vol)$ on $\OGr_{m-2}(M)$ is Hamiltonian.
Indeed, this follows from $[\ham(M,\vol),\ham(M,\vol)]=\ham(M,\vol)$, see
\cite{Li74}, the fact that $\ham(M,\vol)$ acts symplectically 
and the fact that the Lie bracket of two symplectic vector
fields (on $\OGr_{m-2}(M)$)
will be Hamiltonian. In our special situation we do not actually need
this general argument, for we have the following

\begin{lemm}\label{ham}
Let $\mathcal M$ be a connected component of $\OGr_{m-2}(M)$ and choose
$N_0\in\mathcal M$. Then 
$$
\mu:\mathcal M\to\ham(M,\vol)^*,
\quad
\mu(N)(X)=\int_N\alpha-\int_{N_0}\alpha,
\quad
\text{where $i_X\vol=d\alpha$}
$$
is a well defined and injective 
moment map for the $\Hamvol$--action on $\mathcal M$. 
Particularly $\Hamvol$ acts in a Hamiltonian way on $\mathcal M$.
\end{lemm}

\begin{proof}
The definition is meaningful since
$\mu(N)(X)=\int_{\tau}i_X\vol$ 
for any bordism $\tau$ in $M$ with boundary $N-N_0$,
and this expression does not depend on the choice of $\tau$,
for $i_X\nu$ is exact.

The fundamental vector field of $X\in\hamvol$ is
$\zeta_X(N)=X|_N$. To show that $\mu$ is a moment map, we verify that
the function $h(N):=\mu(N)(X)$ is a Hamiltonian function for
the vector field $\zeta_X$. Indeed, up to a constant, $h$ equals 
$\tilde\alpha$, and thus
$$
dh
=d\tilde\alpha
=\widetilde{d\alpha}
=\widetilde{i_X\vol}
=i_{\zeta_X}\tilde\vol
=i_{\zeta_X}\Omega.
$$
The injectivity of this moment map is easily seen choosing
$\alpha$ with appropriate support.
\end{proof}

\begin{prop}\label{trans}
The action of $\Hamvol$ on connected components of 
$\OGr_{n}(M)$ is transitive, provided $m-n\geq2$.
\end{prop}

\begin{proof}
We will show more. Namely we will prove that
$\Hamvol$ acts transitively on every connected component of
$\Emb(N,M)$, the space of embeddings of $N$ in $M$.

First we show that the action of $\Hamvol$ on $\Emb(N,M)$
is infinitesimal transitive,
\ie every vector field along a closed submanifold $N$ in $M$ of codimension 
at least two, can be extended to an exact divergence free vector field on $M$.
We start with an arbitrary extension $Y\in\vf(M)$ of the given vector field
$X\in\Gamma(TM|_N)$. 
By the relative Poincar\'e lemma for the $m$--form $\beta=L_Y\vol$,
there exists an $(m-1)$--form $\lambda$ on 
a tubular neighborhood $U$ of $N$ in $M$,
such that $d\lambda=\beta$ on $U$ and $\lambda|_N=0$.
The relation $i_Z\vol=\lambda$ defines a vector field $Z\in\vf(U)$
with properties: $Z|_N=0$ and $L_Z\vol=\beta$. 
Then $Y-Z$ is a divergence free vector field on $U$ extending $X$.
Since $m-1>n$ we have $H^{m-1}(U)=0$, in particular 
$Y-Z$ is an exact divergence free vector field. 
It can be extended
to an exact divergence free vector field 
$\tilde X\in\vf(M)$ with $\tilde X|_N=X$.

Next we show that every isotopy of $N$ in $M$ extends to an
exact volume preserving diffeotopy of $M$. Indeed,
an isotopy $h_t:N\to M$ determines a smooth family of vector fields $X_t$
on $M$ along $N_t=h_t(N)\subset M$.
By the infinitesimal transitivity
we can extend each $X_t$ to an
exact divergence free vector field $\tilde X_t$.
Looking closer at the construction above, 
we see that the extension can be chosen smoothly depending on $t$.
The diffeotopy $\phi_t$ determined by $\tilde X_t$ is 
exact volume preserving and extends the isotopy $h_t$, \ie $h_t=\phi_t\o h_0$.
So $\Hamvol$ acts transitively on connected components of $\Emb(N,M)$.
\end{proof}

\begin {rema}
The connected components of $\OGr_n(M)$ and $\Emb(N,M)$
can be written as homogeneous spaces of $\Ham(M,\mu)$.
In the first case the isotropy group is the subgroup of
exact volume preserving diffeomorphisms leaving $N\in\OGr_n(M)$ invariant,
in the second case it is the subgroup of exact volume preserving
diffeomorphisms fixing $N$ pointwise.
\end{rema}

Proposition~\ref{kks}, Lemma~\ref{ham}, Proposition~\ref{trans}
and Theorem~\ref{main-theo} prove the following

\begin{theo}\label{coad}
Let $M$ be a closed $m$--dimensional manifold with integral volume form 
$\vol$ and let $\mathcal M$ be a connected component of $\OGr_{m-2}(M)$
equipped with the symplectic form $\Omega=\tilde\vol$. Then there exists
a central extension of $\Hamvol$ by $S^1$ such that $\mathcal M$ is a
coadjoint orbit of this extension. Particularly this coadjoint orbit is
prequantizable.
\end{theo}

\begin{rema}\label{vol_cc}
Recall that the central extension in Theorem~\ref{coad} is the pull-back of 
Kostant's extension by the Hamiltonian action
of $\Ham(M,\nu)$ on $\mathcal M$.
Choose an element $N_0$ in $\mathcal M$.
The moment map $\mu$ from Lemma \ref{ham} vanishes at $N_0$, so by 
$\eqref{cocycle}$ the corresponding Lie algebra cocycle on $\ham(M,\nu)$ is
$c_{N_0}(X,Y)=-\Omega(\zeta_X,\zeta_Y)(N_0)=-\int_{N_0}i_Yi_X\vol$.
The failure of the moment map $\mu$ to be equivariant is
$\kappa(\varphi)(X)=\int_{\varphi(N_0)}\alpha-\int_{N_0}\alpha$,
with $i_X\nu=d\alpha$.
\end{rema}

\begin{rema}
A result of Roger \cite{Ro95} says that 
the second Lie algebra cohomology group of $\hamvol$ 
is isomorphic to $H_{m-2}(M;\R)$,
the 2--cocycle on $\hamvol$ defined by the $(m-2)$--cycle $\sigma$ on $M$ being
$c_\sigma(X,Y)=-\int_\sigma i_Yi_X\vol$.
Every homology class $\sigma$ in $H_{m-2}(M;\Z)$ has a representative 
which is a closed submanifold of codimension 2 in $M$.
The representative $N_0$ can be taken to be the zero set of a section
transversal to the zero section in a rank two vector bundle with Euler class 
the Poincar\'e dual of $\sigma$. 
It follows that all $1$--dimensional central extensions of $\hamvol$
corresponding to $\sigma\in H_{m-2}(M;\Z)$ 
can be integrated to group extensions.
The original construction is due to Ismagilov \cite{I96}.
However, we even get the Lie group structure on the extensions
by using a result in \cite{NV03}. 
\end{rema}

\begin{rema}
Suppose $\mathcal M$ and $\mathcal M'$ are two connected components of
$\OGr_{m-2}(M)$ corresponding to homologous submanifolds of $M$.
Choose $N_0\in\mathcal M$ and $N_0'\in\mathcal M'$. Since $N_0$ and $N_0'$
are homologous we can choose a smooth $(m-1)$--chain $B_0$ in $M$ with
$\partial B_0=N_0'-N_0$. Denote $\goe:=\hamvol$, $G:=\Hamvol$ and define 
$$
\lambda_0\in\goe^*,
\quad\lambda_0(X):=-\int_{B_0}i_X\nu.
$$
This does not depend on the choice of $B_0$, for $i_X\nu$ is exact.
Via Lemma~\ref{ham} $N_0$ and $N_0'$ give rise to moment maps 
$\mu_0:\mathcal M\to\goe^*$ and $\mu_0':\mathcal M'\to\goe^*$
with corresponding cocycles $c_0$ and $c'_0$ and
$\kappa_0:G\to\goe^*$ and $\kappa_0':G\to\goe^*$, 
respectively, see Remark~\ref{vol_cc}.
An easy calculation shows
\begin{equation}\label{iso}
c'_0(X,Y)-c_0(X,Y)=\lambda_0(-[X,Y])
\quad\text{and}\quad
\kappa_0'(\varphi)-\kappa_0(\varphi)=\lambda_0-\varphi\cdot\lambda_0
\end{equation}
for all $X,Y\in\goe$ and all $\varphi\in G$. Again, the unexpected minus sign
stems from the fact, that the usual Lie bracket of vector fields $[X,\cdot]$ 
is an infinitesimal right action, whereas on Lie groups the Lie bracket 
$[X,\cdot]$ is an infinitesimal left action. Moreover the moment maps
give rise to Lie algebra isomorphisms $\R\oplus_{c_0}\goe\to\tilde\goe$ and
$\R\oplus_{c_0'}\goe\to\tilde\goe'$. Using these identifications and
$\lambda_0$ from above we can define a mapping
$$
\Phi_0:\tilde\goe\simeq\R\oplus_{c_0}\goe\to
\R\oplus_{c_0'}\goe\simeq\tilde\goe',
\quad\text{by}\quad
\Phi_0(a,X)=(a+\lambda_0(X),X).
$$
This is an isomorphism of Lie algebras and $G$--equivariant for we have
\eqref{iso}. Particularly the Lie algebra extensions
$0\to\R\to\tilde\goe\to\goe\to0$ and $0\to\R\to\widetilde{\goe'}\to\goe\to0$
are isomorphic, as expected. 

When defining $\Phi_0:\tilde\goe\to\tilde\goe'$ we made two choices, namely
$N_0$ and $N_0'$. We claim that $\Phi_0$ is independent of them. 
Indeed, suppose $N_1\in\mathcal M$
and $N_1'\in\mathcal M'$, choose $B_{01}$ and $B_{01}'$ such that 
$\partial B_{01}=N_1-N_0$ and $\partial B_{01}'=N_1'-N_0'$ and define 
$\rho_{01}\in\goe^*$ by $\rho_{01}(X)=-\int_{B_{01}}i_X\nu$
and $\rho_{01}'\in\goe^*$ by $\rho_{01}'(X)=-\int_{B_{01}'}i_X\nu$,
respectively. Again this does not depend on the
choice of $B_{01}$ or $B_{01}'$. Moreover choose $B_1$, such that
$\partial B_1=N_1'-N_1$ and define $\lambda_1(X):=-\int_{B_1}i_X\nu$.
One easily checks that the composition
$$
\R\oplus_{c_0}\goe\to\tilde\goe\to\R\oplus_{c_1}\goe
\quad\text{is given by}\quad
(a,X)\mapsto(a+\rho_{01}(X),X)
$$
and similarly for $c_0'$, $c_1'$, $\mu_0'$, $\mu_1'$ and $\rho_{01}'$.
Thus $\Phi_0=\Phi_1$ is equivalent to
$\rho_{01}+\lambda_1=\lambda_0+\rho_{01}'$ which is equivalent to
$$
-\int_{B_{01}}i_X\nu
-\int_{B_1}i_X\nu
=-\int_{B_0}i_X\nu
-\int_{B_{01}'}i_X\nu,
$$
but this follows since $i_X\nu$ is exact and the integral is over a cycle.

Summarizing we have seen that whenever the components $\mathcal M$ and
$\mathcal M'$ consist of homologous submanifolds, 
there is a canonic $G$--equivariant
isomorphism of Lie algebras $\Phi:\tilde\goe\to\tilde\goe'$.
Particularly the coadjoint orbits of $\tilde G$ and $\tilde G'$
coincide. 
We are not aware of a more intrinsic definition of $\Phi$, and we don't know
if the corresponding group extensions are isomorphic in this situation.

Finally, suppose $\mathcal M$ is a component of $\OGr_{m-2}(M)$ which 
consists of $0$--homologous submanifolds. Then 
$$
\mu:\mathcal M\to\goe^*,\quad
\mu(N)(X):=\int_N\alpha,\quad\text{with $d\alpha=i_X\nu$}
$$ 
is a $G$--equivariant moment map. So we get a canonic
$G$--equivariant isomorphism of Lie algebras $\tilde\goe\simeq\R\oplus\goe$. 
Moreover $\mathcal M$ is a coadjoint orbit of $G$, canonically.
\end{rema}

\section{Generalized vortex filament equation}

For a Riemannian metric $g$ on $M$ with induced volume form 
$\vol(g)=\vol$, we identify the normal bundle $TN^\perp$ with the 
Riemannian orthonormal bundle $TN^{\perp_g}$ 
and denote by $g$ the induced metric on it.
We endow the symplectic manifold
$(\OGr_{m-2}(M),\Omega)$ with a Riemannian metric
$$
\tilde g( Y _1, Y _2)=\int_Ng( Y _1, Y _2)\vol(g|_N)
\quad\text{for $ Y _1, Y _2\in\Gamma(TN^{\perp_g})$.}
$$
For $N\in\OGr_{m-2}(M)$, the vector bundle $TN^{\perp_g}$
is oriented, 2--dimensional and has a metric, 
so we can define a fiber wise complex structure $J$ on $TN^{\perp_g}$ by 
rotation with $+90$ degrees. It induces an almost complex structure $\tilde
J$ on 
$\OGr_{m-2}(M)$ which is compatible with $\Omega$ and $\tilde g$,
that is $\Omega( Y _1, Y _2)=\tilde g(\tilde J Y _1, Y _2)$.

The $g$--volume of the submanifold gives a smooth function on $\OGr_{m-2}(M)$
\begin{equation}\label{hamfunc}
h:\OGr_{m-2}(M)\to\R,\quad
h(N)=\int_N\vol(g|_N).
\end{equation}

\begin{lemm}\label{vortex_prop}
For the $\tilde g$--gradient of $h$ we have $(\grad h)(N)=-\tr\II_N$, 
where $\II_N\in\Gamma(S^2T^*N\otimes TN^{\perp_g})$ 
denotes the second fundamental form of the submanifold $N$.
\end{lemm}

\begin{proof}
For $ Y \in\Gamma(TN^{\perp_g})$ we have
\begin{align*}
dh( Y )
&=\frac 12\int_N\tr(L_ Y  g)\vol(g|_N)
\\
&=\int_N\tr(\nabla Y )\vol(g|_N)
\\
&=-\int_N\tr g(\II_N, Y )\vol(g|_N)
\\
&=-\tilde g(\tr\II_N, Y ).
\end{align*}
Since $\tilde g$ is weakly non-degenerated we conclude $\grad h=-\tr\II$.
\end{proof}

Since $\tilde J$, $\Omega$ and $\tilde g$ are compatible,
the Hamiltonian vector field of $h$ is
$X_h=\tilde J(\grad H)=J\tr\II$ and this proves the following

\begin{prop}
The Hamiltonian equation for the Hamiltonian function \eqref{hamfunc} is
$$
\frac\partial{\partial t}N_t=J\tr\II_{N_t}.
$$
\end{prop}

In dimension $m=3$ this equation is known as the vorticity filament
equation, see \cite{MW83}. 

\begin{rema}
Let $N$ be a closed oriented manifold of dimension $m-2$.
The expression $J\tr\II$ can also be considered as a vector field on 
$\Emb(N,M)$. Suppose $\iota_t$ is a curve of embeddings solving 
\begin{equation}
\label{vortex_equ}
\frac\partial{\partial t}\iota_t=J\tr\II_{\iota_t(N)}
\end{equation} 
and let $f$ be an orientation preserving diffeomorphism
of $N$. Then $\iota_t\circ f$ will again be a solution of \eqref{vortex_equ}.
The geometric interpretation of this fact is the following. When restricting
to suitable connected components, the space of embeddings becomes a principal
$\Diff(N)$--bundle over the non-linear Grassmannian. Here $\Diff(N)$ denotes
the group of orientation preserving diffeomorphisms. Using the Riemannian
metric we can write down a connection of this bundle,
known as a mechanical connection. For an embedding
$\iota:N\to M$, the vertical tangent space is the space of vector fields
along $\iota$ tangent to $\iota(N)$. So the space of vector fields along
$\iota$ having values in the Riemannian orthogonal complement of $\iota(N)$
is a complement to the vertical tangent space. This complements define a 
connection, which is obviously a principal connection. Regarding the expression
$J\tr\II$ as a vector field on the space of embeddings, just means 
considering the horizontal lift of $J\tr\II$. 
Since the connection is principal, parallel transport will be
$\Diff(N)$--equivariant. This translates to $\iota_t\circ f$ is a solution of
\eqref{vortex_equ} iff $\iota_t$ was. 
\end{rema}

\begin{rema}
Let $\iota_t$ be a curve of embeddings in $M$. Then we get a curve of
Riemannian metrics $\iota_t^*g$ on $N$. This gives rise to a curve of volume
forms $\vol(\iota_t^*g)$ on $N$. If $\iota_t$ is a solution of
\eqref{vortex_equ} this curve will be constant. Indeed, for every horizontal
$\iota_t$ one shows $\frac\partial{\partial t}\vol(\iota_t^*g)
=-g(\tr\II,\frac\partial{\partial t}\iota_t)$
as in the proof of Lemma~\ref{vortex_prop}.
If $\iota_t$ solves \eqref{vortex_equ} this implies
$\frac\partial{\partial t}\vol(\iota_t^*g)
=-g(\tr\II,\frac\partial{\partial t}\iota_t)
=g(J\frac\partial{\partial t}\iota_t,\frac\partial{\partial t}\iota_t)=0.$
In the case of oriented knots in a $3$--dimensional $M$,
this implies that a solution of \eqref{vortex_equ},
parameterized by arc length at time $t_0$, will 
have the same property for every time $t$.
\end{rema}

\section{Symplectic Grassmannians}

Suppose $(M,\omega)$ is a closed connected symplectic manifold.
Let $\SOGr_{2k}(M)\subseteq\OGr_{2k}(M)$ denote the open
subset of oriented submanifolds which are symplectic. 
We don't assume the elements in $\SOGr_{2k}(M)$ to be oriented by
their symplectic form.
Note that $\SOGr_{2k}(M)$ is invariant under the action of
the group of symplectic diffeomorphisms $\Diff(M,\omega)$.
Set $\alpha:=\omega^{k+1}$. 
Then $\Omega:=\tilde\alpha$ is a closed $2$--form on
$\OGr_{2k}(M)$. Note that $\SOGr_{2k}(M)\subseteq\OGr_{2k}(M)$ is an open
subset on which $\Omega$ is (weakly) non-degenerate,
hence a symplectic manifold.
Indeed, for an almost complex structure $J$ on $M$ tamed by $\omega$
and $Y\in\Gamma(TN^{\perp_\omega})$ we have
$\Omega_N(Y,JY)=(k+1)\int_N\omega(Y,JY)\omega^k$, vanishing iff $Y=0$.

Let $\Ham(M,\omega)$ denote the Lie group of Hamiltonian diffeomorphisms
with Lie algebra $\ham(M,\omega)$ of Hamiltonian vector fields on $M$.
The action of $\Ham(M,\omega)$ on $\SOGr_{2k}(M)$
is Hamiltonian. Indeed we already know that the action is symplectic and
since $[\ham(M,\omega),\ham(M,\omega)]=\ham(M,\omega)$, see \cite{C70},
the action must be Hamiltonian. In our special situation one does not have
to use this general argument, for one can write down Hamilton functions.

\begin{lemm}\label{sympinj}
The mapping
$$
\mu:\SOGr_{2k}(M)\to\ham(M,\omega)^*,
\quad
\mu(N)(X):=(k+1)\int_Nf\omega^k, 
$$
is an injective equivariant moment map for the 
$\Ham(M,\omega)$--action on $\SOGr_{2k}(M)$. 
Here $f$ is the unique Hamilton function of $X$ with zero integral.
Particularly the action is Hamiltonian.
\end{lemm}

\begin{proof}
First we show that $h(N):=(k+1)\int_Nf\omega^k$ is 
a Hamiltonian function for the fundamental
vector field $\zeta_X$ of $X\in\ham(M,\omega)$. Note that
$h=(k+1)\widetilde{f\omega^k}$ and 
$(k+1)d(f\omega^k)
=(k+1)df\wedge\omega^k
=i_X\omega^{k+1}$. Thus
$$
dh
=(k+1)d\widetilde{f\omega^k}
=(k+1)\widetilde{d(f\omega^k)}
=\widetilde{i_X\omega^{k+1}}
=i_{\zeta_X}\widetilde{\omega^{k+1}}
=i_{\zeta_X}\Omega.
$$
So $\mu$ is a moment map. The injectivity is obvious.
Finally for every $\varphi\in\Ham(M,\omega)$ we have
$$
\mu(\varphi(N))(X)=(k+1)\int_{\varphi(N)}f\omega^k=
(k+1)\int_N(\varphi^*f)\omega^k=\mu(N)(\varphi^*X).
$$
and thus $\mu$ is equivariant.
\end{proof}

\begin{prop}\label{symptrans}
The group $\Ham(M,\omega)$ acts transitively on every connected component 
of $\SOGr_{2k}(M)$.
\end{prop}

\begin{proof}
We first show that the action is infinitesimal transitive. So suppose
$N\in\SOGr_{2k}(M)$ and let $X$ be a tangent vector at $N$, \ie a section of
the normal bundle $TN^\perp$.
Since $N$ is a symplectic submanifold we can identify the normal
bundle with the $\omega$--orthogonal complement
$TN^{\perp_\omega}$ of $TN$. So we may assume that $X$ is a section 
of $TN^{\perp_\omega}$. Consider $i_X\omega$ as a function, say $\lambda$, on 
the total space $E$ of $TN^{\perp_\omega}$ which happens to be linear 
along the fibers. One easily shows that
$d\lambda=i_X\omega$ along $N\subseteq E$. Considering $E$
as a tubular neighborhood of $N$ one easily gets a
function $\lambda'$ on $M$ such that
$d\lambda'=i_X\omega$ along $N\subseteq M$. So the Hamiltonian vector 
field to $\lambda'$ will be an extension of $X$. Thus
$\Ham(M,\omega)$ acts infinitesimally transitive on $\SOGr_{2k}(M)$.

Suppose $N_t$ is a curve in $\SOGr_{2k}(M)$ and set 
$X_t:=\frac\partial{\partial t}N_t$, a section of $TN_t^\perp$.
For every fixed time $t$ the section $X_t$ can be extended to a vector
field in $\ham(M,\omega)$ as shown above. Moreover it is clear that this
extension can be chosen smoothly with respect to the parameter $t$. Now the
flow of this extension clearly gives a curve in $\Ham(M,\omega)$
transporting, say, $N_0$ to $N_1$.
\end{proof}

Proposition~\ref{kks}, Lemma~\ref{sympinj}, Proposition~\ref{symptrans} 
and Theorem~\ref{main-theo} prove the following

\begin{theo}
Let $(M,\omega)$ be a symplectic manifold, such that $[\omega]^{k+1}\in
H^{2k+2}(M;\R)$ is
integral and let $\mathcal M$ denote a connected component of $\SOGr_{2k}(M)$
endowed with the symplectic form $\Omega=\widetilde{\omega^{k+1}}$.
Then $\mathcal M$ is a coadjoint orbit of $\Ham(M,\omega)$.
Particularly this coadjoint orbit is prequantizable.
\end{theo}

\begin{rema}
Since we have an equivariant moment map the Lie algebra extension
$0\to\R\to\tilde\goe\to\goe\to0$ from section~\ref{univ_constr} 
with $\goe=\ham(M,\omega)$ is trivial. 
This is the reason why $\SOGr_{2k}(M)$ can be considered as coadjoint orbit 
of $\Ham(M,\omega)$ rather than as coadjoint orbit of a central
extension. However, the group extension 
\begin{equation}\label{sgg}
1\to S^1\to\tilde G\to G\to1
\end{equation}
from section~\ref{univ_constr} with $G=\Ham(M,\omega)$
may very well be non-trivial as the following 
example shows.
\end{rema}

\begin{exam}\label{gru_ext_ex}
Let $M=S^2$, $\omega$ the standard symplectic form of mass $1$ and let
$k=0$. Then $\mathcal M=S^2$ is a connected component of $\SOGr_{2k}(M)$
and $\Omega=\omega$, \cf Example~\ref{trivial_ex}.
The bundle $\mathcal P\to\mathcal M$ is the Hopf fibration $S^3\to S^2$
and $\eta$ is the standard contact structure on $S^3$.
Since $M=\mathcal M$, the group extension \eqref{sgg}
is trivial iff Kostant's extension 
$$
1\to S^1\to\Aut(S^3,\eta)\to\Ham(S^2,\omega)\to 1
$$ 
is trivial. The equivariant moment map from Lemma~\ref{sympinj} 
provides a Lie algebra homomorphism
$\sigma:\ham(S^2,\omega)\to\aut(S^3,\eta)$, right inverse to the projection.
So Kostant's group extension is trivial iff this Lie algebra homomorphism
integrates to a group homomorphism. However this is not the case. Indeed,
the loop in $\Ham(S^2,\omega)$ given by rotation around an axis
does not integrate to a closed curve in $\Aut(S^3,\eta)$.
To see this, note first that the Hamilton function generating the rotation
vanishes along the equator, for it has zero integral. So $\sigma$ maps the
Hamilton vector field to an element of $\aut(S^3,\eta)$ which is horizontal
over the equator of $S^2$. So integrating our loop of rotation gives a 
curve in $\Aut(S^3,\eta)$ whose flow lines over the equator of $S^2$ are
horizontal. Such a flow line has holonomy $1/2$, for this
is the total curvature of a hemisphere. Thus it is not closed.

Alternatively one can use the fact that the Hamilton function generating the
rotation has values $\pm1/2$ at the poles.
\end{exam}


\begin{thebibliography}{MW83}

\bibitem[C70]{C70}
E. Calabi,
{\it On the group of automorphisms of a symplectic manifold},
Problems in Analysis, Symp. in honor of S. Bochner, 1--26,
Princeton University Press, 1970.

\bibitem[I96]{I96}
R. S. Ismagilov,
{\it Representations of infinite-dimensional groups},
Translations of Mathematical Monographs {\bf 152}, 
American Mathematical Society, Providence, RI, 1996.

\bibitem[K70]{Ko70}
B. Kostant,
{\it Quantization and unitary representations},
Lectures in modern analysis and applications III, 87--208,
Lecture Notes in Math. {\bf 170}, Springer, Berlin, 1970.

\bibitem[KM97]{KM97}
A. Kriegl and P. W. Michor,
{\it The convenient setting of global analysis},
Mathematical Surveys and Monographs {\bf 53},
American Mathematical Society, Providence, RI, 1997.

\bibitem[L74]{Li74}
A. Lichnerowicz, 
{\it Alg\`ebre de Lie des automorphismes 
infinit\'esimaux d'une structure unimodulaire},
Ann. Inst. Fourier {\bf 24}(1974), 219--266.

\bibitem[MW83]{MW83}
J. Marsden and A. Weinstein,
{\it Coadjoint orbits, vortices, and Clebsch variables for incompressible
fluids},
Phys. D {\bf 7}(1983), 305--323.

\bibitem[NV]{NV03}
K.--H. Neeb and C. Vizman,
{\it Flux homomorphisms and principal bundles
over infinite dimensional manifolds}, 
to appear in Monatsh. Math.

\bibitem[R95]{Ro95}
C. Roger, {\it Extensions centrales d'alg\`ebres et de groupes de Lie
de dimension infinie, alg\`ebre de Virasoro et g\'en\'eralisations},
Rep. Math. Phys. {\bf 35}(1995), 225--266.

\end{thebibliography}
\end{document}